\documentclass[12pt]{article}

\usepackage{amsmath,amssymb,amsthm}
\usepackage{amsfonts}
\usepackage[mathscr]{eucal}
\begin{document}
\newtheorem{thm}{Theorem}
\numberwithin{thm}{section}
\newtheorem{lemma}[thm]{Lemma}
\newtheorem{remark}{Remark}
\newtheorem{corr}[thm]{Corollary}
\newtheorem{proposition}{Proposition}
\newtheorem{theorem}{Theorem}[section]
\newtheorem{deff}[thm]{Definition}
\newtheorem{case}[thm]{Case}
\newtheorem{prop}[thm]{Proposition}
\numberwithin{equation}{section}
\numberwithin{remark}{section}
\numberwithin{proposition}{section}
\newtheorem{corollary}{Corollary}[section]
\newtheorem{others}{Theorem}
\newtheorem{conjecture}{Conjecture}\newtheorem{definition}{Definition}[section]
\newtheorem{cl}{Claim}
\newtheorem{cor}{Corollary}
\newcommand{\ds}{\displaystyle}

\newcommand{\stk}[2]{\stackrel{#1}{#2}}
\newcommand{\dwn}[1]{{\scriptstyle #1}\downarrow}
\newcommand{\upa}[1]{{\scriptstyle #1}\uparrow}
\newcommand{\nea}[1]{{\scriptstyle #1}\nearrow}
\newcommand{\sea}[1]{\searrow {\scriptstyle #1}}
\newcommand{\csti}[3]{(#1+1) (#2)^{1/ (#1+1)} (#1)^{- #1
 / (#1+1)} (#3)^{ #1 / (#1 +1)}}
\newcommand{\RR}[1]{\mathbb{#1}}
\thispagestyle{empty}
\begin{titlepage}
\title{\bf Iterated Brownian motion in bounded domains in $\RR{R}^{n}$}
\author{Erkan Nane\thanks{ Supported in part by NSF Grant \# 9700585-DMS }\\
Department of Mathematics\\
Purdue University\\
West Lafayette, IN 47906 \\
enane@math.purdue.edu}
\maketitle
\begin{abstract}
\noindent {\it Let  $\tau _{D}(Z) $ is the first
exit time of iterated Brownian motion from a domain $D \subset
\RR{R}^{n}$ started at $z\in D$ and let $P_{z}[\tau _{D}(Z)  >t]$ be
its distribution.   In this paper
 we establish the exact asymptotics of $P_{z}[\tau _{D}(Z)  >t]$
 over bounded domains as an extension of the result in DeBlassie
 \cite{deblassie}, for $z\in D$
$$
P_{z}[\tau_{D}(Z)>t]\approx t^{1/2}
\exp(-\frac{3}{2}\pi^{2/3}\lambda_{D}^{2/3}t^{1/3}) ,\ as\ t\to\infty .
$$
 We also study  asymptotics of the life time of Brownian-time Brownian
 motion (BTBM),
$Z^{1}_{t}=z+X(Y(t))$, where $X_{t}$ and $Y_{t}$ are independent
one-dimensional Brownian motions.}

\end{abstract}

\end{titlepage}
\section{Introduction and statement of main results}

Properties of iterated Brownian motion (IBM) analogous to the
properties of Brownian motion have been studied extensively by
several authors \cite{allouba2,  allouba1, bandeb,
burdzy1, burdzy2, bukh, csaki, deblassie, eisenbaumshi, klewis,
 nane, xiao}.  Several other iterated processes including
 Brownian-time Brownian motion (BTBM) have also been studied \cite{allouba2, allouba1,
 koslew}. One of the main differences between these iterated processes  and Brownian motion is
that they  are not Markov processes.  However, these
processes have connections
with the
 parabolic operator
$\frac{1}{8} \Delta ^{2} -\frac{\partial}{\partial t}$, as
described in \cite{allouba1,deblassie}.


To define the iterated  Brownian motion $Z_{t}$
started at $z \in \RR{R}$,
 let $X_{t}^{+}$, $X_{t}^{-}$ and $Y_{t}$
 be three  independent
one-dimensional Brownian motions, all started at $0$. Two-sided
Brownian motion is defined by
\[ X_{t}=\left\{ \begin{array}{ll}
X_{t}^{+}, &t\geq 0\\
X_{(-t)}^{-}, &t<0.
\end{array}
\right. \] Then the iterated Brownian motion started at $z \in
\RR{R}$ is
\[ Z_{t}=z+X(Y_{t}),\ \ \    t\geq 0.\]

In $\RR{R}^{n}$, one requires $X^{\pm}$ to be independent
$n-$dimensional Brownian motions. This is the version of the
iterated Brownian motion due to Burdzy, see \cite{burdzy1}.

In what follows  we will write $f \approx g$ and $f \lesssim g$ to
mean that for some positive $C_{1}$ and $C_{2}$, $C_{1}\leq f/g
\leq C_{2}$ and $f \leq C_{1} g$, respectively. We will also write
$f(t) \sim g(t)$, as $t\rightarrow \infty $,  to mean that $f(t) / g(t)
\rightarrow 1$, as $t\rightarrow \infty $.

Let $\tau_{D}$ be the first exit time of Brownian motion from a
domain $D\subset \RR{R}^{n} $. The large time behavior of
$P_{z}[\tau_{D}>t]$ has been studied for  several types of domains,
 including general cones
\cite{bansmits, deblassie1}, parabola-shaped domains \cite{bds,
lshi},
 twisted domains \cite{DSmits} and
bounded domains \cite{blab}. Our aim in this article is to do the same for the exit time of IBM over bounded domains in $\RR{R}^{n}$ and for the exit times of BTBM over several domains in  $\RR{R}^{n}$.

In particular, the large time asymptotics of the lifetime of
Brownian motion in general cones has been studied by several
people including  Burkholder \cite{burkholder}, DeBlassie \cite{deblassie1} and  Ba\~{n}uelos and
Smits \cite{bansmits}. Let $D$ be an open cone with vertex $0$
such that $S^{n-1}\cap D$ is regular for the Laplace-Beltrami
operator $L_{S^{n-1}}$ on the sphere $S^{n-1}$. Then for some $p(D)>0$ (see
\cite{deblassie1} and \cite{bansmits})
$$
P_{x}[\tau_{D}>t]\sim C(x)t^{-p(D)},\ as \ t\to\infty.
$$
Now let $D\subset \RR{R}^{n}$. Let $ \tau_{D}(Z)=\inf \{ t\geq 0:\
Z_{t}\notin D\},$ be the first exit time of $Z_{t}$ from $D$. When
$D$ is a generalized cone, using the results of Ba\~{n}uelos and
Smits, DeBlassie \cite{deblassie} obtained; for $z\in D$, as
$t\to\infty$,
$$
P_{z}[\tau_{D}(Z)>t]\approx \left\{ \begin{array}{ll}
t^{-p(D)}, &p(D)<1\\
t^{-1}\ln t, &p(D)=1\\
t^{-(p(D)+1)/2} ,& p(D)>1.
\end{array}
\right.
$$
For  parabola-shaped domains the study of exit time asymptotics
for Brownian motion was initiated  by Ba\~{n}uelos, DeBlassie and
Smits \cite{bds} to answer the question: Are there domains in
$\RR{R}^n$ for which the distribution of the exit time is
sub-exponential?
 They showed  that for  the
 parabola $\mathcal{P}=\{(x,y):
x>0, |y|<A\sqrt{x} \}$, $A>0$  there exist positive constants
$A_{1}$ and $A_{2}$ such that for $z\in \mathcal{P}$,
\begin{eqnarray}
-A_{1} & \leq & \liminf_{t \rightarrow \infty} \ t^{-\frac{1}{3}}
\log \ P_{z}[\tau _{\mathcal{P}} >t]
 \leq
  \limsup_{t \rightarrow \infty}\  t^{-\frac{1}{3}} \log
 P_{z}[\tau _{\mathcal{P}} >t] \leq  -A_{2}. \nonumber
\end{eqnarray}
Subsequently,   Lifshits and Shi \cite{lshi} found that the above
limit exists for parabola-shaped domains $P_{\alpha}=\{ (x,Y)\in
\RR{R} \times \RR{R}^{n-1}: x>0, |Y|<Ax^{\alpha} \}$, $0 <\alpha
<1$ and $A>0$ in any dimension; for $z\in P_{\alpha}$,
\begin{equation}
\lim_{t\to\infty}t^{-(\frac{1-\alpha } { 1+\alpha} )}
\log  P_{z}[\tau _{\alpha}> t]=- l,
\end{equation}
where
\begin{equation}\label{brownianlimit}
l=(\frac{1+\alpha}{ \alpha})\left( \frac{\pi \jmath_{(n-3)/2}^{2/ \alpha}}{A^{2}2^{(3\alpha+1)/ \alpha}
((1-\alpha)/ \alpha)^{(1-\alpha)/ \alpha }}
\frac{\Gamma^{2}(\frac{1- \alpha }{2\alpha})}{\Gamma^{2}(\frac{1}{2\alpha})}
\right)^{\frac{\alpha }{(\alpha +1)}}.
\end{equation}
Here  $\jmath _{(n-3)/ 2}$ denotes the smallest positive zero of the Bessel function
$J_{(n-3)/2}$ and
$\Gamma $ is the Gamma function.

Using the results for Brownian motion in parabola-shaped domains
we established in \cite{nane} with $l$  given by (\ref{brownianlimit}), for $z \in P_{\alpha}$,
\begin{eqnarray}
  \lim_{t\to\infty}t^{-(\frac{ 1- \alpha }{3+\alpha})}\log P_{z}[\tau _{\alpha}(Z)> t]
 = -(\frac{3+\alpha }{2+ 2 \alpha } )
(\frac{1+\alpha }{1-\alpha})^{(\frac{1- \alpha }{3+ \alpha} )} \pi
^{(\frac{2-2\alpha }{3+\alpha } )} l^{  (\frac{2+2\alpha
}{3+\alpha }) }.  \nonumber
\end{eqnarray}

For many bounded domains  $D\subset \RR{R}^{n} $ the
asymptotics of $P_{z}[\tau_{D}>t]$ is well-known (See \cite{blab} for a more precise statement of this.) For $z\in D$,
\begin{equation}\label{exp}
\lim_{t\to \infty} e^{\lambda_{D} t}P_{z}[\tau_{D}>t]=
\psi(z)\int_D \psi(y) dy,
\end{equation}
where $\lambda_{D}$ is the first eigenvalue of $\frac{1}{2} \Delta$ with Dirichlet
boundary conditions  and $\psi$ is its corresponding
eigenfunction.

In \cite{deblassie}, DeBlassie proved  in the case of iterated
Brownian motion in bounded domains for $z\in D$,
\begin{equation}\label{boundeddomain0}
\lim_{t\to \infty}t^{-1/3}\log P_{z}[\tau_{D}(Z)>t]=\ -\frac{3}{2}
\pi^{2/3} \lambda_{D}^{2/3}.
\end{equation}

The limits (\ref{exp}) and (\ref{boundeddomain0}) are very different in that the latter
involves taking the logarithm which may kill many unwanted terms in the exponential.
It is then natural to ask if it is possible to obtain an analogue of (\ref{exp}) for
IBM. That is, to remove the log in (\ref{boundeddomain0}).  In this paper we prove the
following theorem.
\begin{theorem}\label{boundeddomain1}
Let $D\subset \RR{R}^{n}$ be bounded domain for which (\ref{exp})
holds pointwise and let $\lambda_D$ and $\psi$ be as above. Then
for $z\in D$,
\begin{eqnarray}
 2C(z)  & \leq &
\liminf_{t\to\infty} t^{-1/2}
\exp(\frac{3}{2}\pi^{2/3}\lambda_{D}^{2/3}t^{1/3})P_{z}[\tau_{D}(Z)>t]\nonumber\\
& \leq  & \limsup_{t\to\infty}
t^{-1/2}\exp(\frac{3}{2}\pi^{2/3}\lambda_{D}^{2/3}t^{1/3})
P_{z}[\tau_{D}(Z)>t]\leq \pi C(z),\nonumber
\end{eqnarray}
where $ C(z)=\lambda_{D}\sqrt{2\pi /3}\left(
\psi(z)\int_{D}\psi(y)dy\right) ^{2}. $
\end{theorem}
We also obtain a version of the above Theorem \ref{boundeddomain1}
 for another closely
related process, the so called
 Brownian-time Brownian motion (BTBM).  To define this,
 let $X_{t}$ and $Y_{t}$
 be two  independent
one-dimensional Brownian motions, all started at $0$. BTBM is
defined to be $Z_{t}^{1}=x+X(|Y_{t}|)$. Properties of this process
and its connections to PDE's have been studied in
\cite{allouba2}, \cite{allouba1} and  \cite{koslew}.
 Analogous to Theorem \ref{boundeddomain1} we have
the following result for this process.
\begin{theorem}\label{boundeddomain2}
Let $D\subset \RR{R}^{n}$,
 $\lambda_{D}$ and $\psi$ be as in the statement of Theorem \ref{boundeddomain1}.
 Let $\tau_{D}(Z^{1})$ be the first exit time of BTBM from $D$. Then for $z\in D$,
$$\lim_{t\to\infty} t^{-1/6}\exp (\frac{3}{2} 2^{-2/3}
\pi^{2/3} \lambda_{D}^{2/3}t^{1/3})P_{z}[\tau_{D}(Z^{1})>t]
=C(\lambda_{D})\psi (z)\int_{D}\psi (y)dy ,$$ where
$C(\lambda_{D})=\pi^{-1/6}2^{13/6}3^{-1/2}\lambda_{D}^{1/3}.$ This
limit is uniform on compact subsets of $D$.
\end{theorem}
Notice that the limits in Theorems \ref{boundeddomain1} and \ref{boundeddomain2}  are different, even at the exponential level.

We obtain the following inequality between distributions of
$\tau_{D}(Z)$ and $\tau_{D}(Z^{1})$.
\begin{theorem}\label{IBM1IBM2}
Let $D\subset \RR{R}^{n}$. Then for all $z\in D$ and all $t>0$,
$$
P_{z}[\tau_{D}(Z)>t]\ \leq \ 2 P_{0}[\tau_{D}(Z^{1})>t].
$$
\end{theorem}
\begin{remark}
Notice that, from the theorems  proved in this paper, the reverse
inequality in Theorem \ref{IBM1IBM2} cannot hold for all large
$t$, in the case of domains $D \subset \RR{R}^{n}$ considered
(i.e. bounded domains with regular boundary, parabola-shaped
domains, twisted domains.)
\end{remark}

The paper is organized as follows. In \S 2 we give some
preliminary lemmas to be used in the proof of main results.
Theorem \ref{boundeddomain1} is proved in \S3. \S4 is devoted to
prove Theorem \ref{boundeddomain2} and some other results on the
exit time asymptotics of BTBM over several domains. In \S5, we
compare the exit time distributions of IBM and BTBM. In \S6, we
prove several asymptotic results to be used in the proof of main
results.
\section{Preliminaries}
In this section we state some preliminary facts that will be
used in the proof of main results.

 The main fact is the following Tauberian theorem (\cite[Laplace transform method, 1958, Chapter 4]{debruijn}). Laporte \cite{laporte} also studied this type 
of integrals. 
Let $h$  and $f$ be continuous functions on $\RR{R}$. Suppose $f$
is non-positive and has a global max at $x_{0}$, $f'(x_{0})=0$, $
f''(x_{0})<0$ and $h(x_{0})\neq 0$  and $\int_{-\infty}^{\infty} h(x)\exp(\lambda f(x)) <\infty$ for all $\lambda>0$. Then as $\lambda \to\infty$,
\begin{equation}\label{laplacemethod}
\int_{0}^{\infty}h(x)\exp(\lambda f(x))dx \sim
h(x_{0})\exp(\lambda f(x_{0}))\sqrt{\frac{2\pi}{\lambda
|f''(x_{0})|}}.
\end{equation}

It can be easily seen from Laplace transform method that as $\lambda \to \infty$,
\begin{equation}\label{generalasymptotic}
\int_{0}^{\infty}\exp (-\lambda (x+x^{-2}))dx\sim \exp(-3 \lambda
2^{-2/3})\sqrt{\frac{2^{4/3}\pi}{3 \lambda}}.
\end{equation}

Similarly, as $t \to \infty$,
\begin{equation}\label{newasymptotic}
\int_{0}^{\infty}\exp (-\frac{at}{u^{2}}-bu)du\ \sim
\sqrt{\frac{\pi}{3}}2^{2/3}a^{1/6}b^{-2/3}t^{1/6} \exp (-3
a^{1/3}b^{2/3}2^{-2/3}t^{1/3}).
\end{equation}

This  follows from equation
(\ref{generalasymptotic}) and after making the change of variables
$u=(atb^{-1})^{1/3}x$.

Finally, we obtain, as $t \to \infty$,
\begin{equation}\label{generalast11}
\int_{0}^{\infty}u\exp (-\frac{at}{u^{2}}-bu)du\ \sim 2\sqrt{\frac{\pi}{
3}}a^{1/2} b^{-1} t^{1/2}\exp (-3 a^{1/3}b^{2/3}2^{-2/3}t^{1/3}).
\end{equation}

\section{Iterated Brownian motion in bounded domains}
If $D\subset \RR{R} ^{n}$  is an open set, write
\[
\tau_{D}^{\pm}(z)=\inf \{ t\geq 0:\ \ X_{t}^{\pm} +z \notin D\},\]
and if $I\subset \RR{R}$ is an open interval, write
\[
\eta _{I}=\eta (I)= \inf \{ t\geq 0 :\ \  Y_{t}\notin I\}.
\]
Recall that $\tau _{D}(Z)$ stands for the first exit time of
iterated Brownian motion from  $D$. As in DeBlassie
\cite[\S3.]{deblassie}, we have by the continuity of the paths for
$Z_{t}=z+X(Y_{t})$, if $f$ is the probability density of $\tau
_{D}^{\pm}(z)$  
\begin{equation}\label{translation}
P_{z}[\tau _{D}(Z) > t]=\int_{0}^{\infty} \!\int_{0}^{\infty}
 P_{0}[\eta _{(-u,v)}>t]  f(u) f(v)dvdu.
\end{equation}

\begin{proof}[The proof of Theorem \ref{boundeddomain1}]
The following is well-known 
\begin{equation}\label{distr1}
P_{0}[\eta _{(-u,v)} >t]=\frac{4}{\pi
}\sum_{n=0}^{\infty}\frac{1}{2n+1} \exp (-\frac{(2n+1)^{2}\pi
^{2}}{2(u+v)^{2}} t) \sin \frac{(2n+1)\pi u}{u+v}, 
\end{equation}
(see Feller \cite[pp. 340-342]{feller}).\newline Let $\epsilon
>0$. From Lemma \ref{lemmaA.1}, choose $M>0$ so large that
\begin{equation}\label{papprox1}
 (1-\epsilon)  \frac{4}{\pi} e ^{- \frac{\pi ^{2} t}{2} } \sin \pi x \leq P_{x}[\eta _{(0,1)} >t] 
 \leq (1+\epsilon)  \frac{4}{\pi} e ^{- \frac{\pi ^{2} t}{2} } \sin \pi x,  
\end{equation}
 for $t\geq M $, uniformly $ x \in (0,1)$.
Let $0<\delta <1/2 $, from the Jordan inequality for the sine
function in the interval $(0,\pi/2]$,
\begin{equation}
2x \leq \sin \pi x \leq \pi x, \ \   x\in
(0,\delta].\label{sinapprox1}
\end{equation}
For a bounded domain with regular boundary it is well-known (see
\cite[page 121-127]{blab}) that there exists an increasing sequence of
eigenvalues, $\lambda_{1}< \lambda_{2} \leq \lambda_{3}\cdots ,$
and eigenfunctions $\psi_{k}$ corresponding to $\lambda_{k}$ such
that,
\begin{equation}\label{eigenvalue-ex0}
P_{z}[\tau_{D}\leq t]=\sum_{k=1}^{\infty}
\exp (-\lambda_{k} t)\psi_{k}(z)\int_{D}\psi_{k}(y)dy.
\end{equation}
From the arguments in DeBlassie \cite[Lemma A.4]{deblassie}
\begin{equation}\label{eigenvalue-ex}
f(t)=\frac{d}{dt}P_{z}[\tau_{D}\leq t]=\sum_{k=1}^{\infty}
\lambda_{k} \exp (-\lambda_{k} t)\psi_{k}(z)\int_{D}\psi_{k}(y)dy.
\end{equation}
 Finally choose $K>0$ so large that
 $$A(z)(1-\epsilon)\exp (-\lambda_{D} u)\leq f(u)\leq A(z)(1+\epsilon)\exp (-\lambda_{D} u) $$
for all $u\geq K$, where
$$A(z)=\lambda_{1}\psi_{1}(z)\int_{D}\psi_{1}(y)dy=\lambda_{D}\psi (z)\int_{D}\psi (y)dy.$$
We further assume that $t$ is so large that $K<\delta \sqrt{t/M}$.
 Define $A$ for $\delta < 1 /2 ,\  K >0 $ and $ M>0$ as
\[
A=\left\{(u,v):\ \ K\leq u\leq \delta \sqrt{\frac{t}{M}} ,\
\frac{1-\delta}{\delta} u \leq v\leq \sqrt{\frac{t}{M}}-u
\right\}.
\]
On the set $A$, since $\delta < 1/2 $, we have $v\geq
(\frac{1}{\delta} -1)u>u>K$ and $u+v >\frac{u}{\delta}$, this
gives $\frac{u}{u+v}\leq \delta. $

 By equations  (\ref{papprox1}) and (\ref{sinapprox1}), $P_{z}[\tau _{D}(Z)>t]=P[\eta _{(-\tau _{D}^{-}(z), \tau _{D}^{+}(z))}>t]$ is
\begin{eqnarray}
\geq C^{1} \int_{K}^{\delta \sqrt{t/ M}} \int_{(1-\delta)u/
\delta}^{\sqrt{t/ M}-u} \frac{u }{ (u+v)}\exp (-\frac{\pi ^{2}
t}{2 (u+v)^{2}})
 \exp (-\lambda_{D} (u+v)) dvdu, \nonumber
\end{eqnarray}
where $C^{1}=C^{1}(z)=4(4/\pi)A(z)^{2}(1-\epsilon)^{3}$. Changing the
variables $x=u+v, z=u$ the integral is
\[
= C^{1}\int_{K}^{\delta \sqrt{t/ M}} \int_{z/ \delta}^{\sqrt{t/
M}}\frac{z}{x} \exp(-\frac{\pi ^{2} t}{2 x^{2}})  \exp
(-\lambda_{D} x) dxdz,
\]
and reversing the order of integration
\begin{eqnarray}
\ &\ & =C^{1}\int_{K/ \delta}^{\sqrt{t/ M}} \int_{K}^{\delta x}\
\frac{z}{x}\exp(-\frac{\pi ^{2} t}{2 x^{2}})  \exp (-\lambda_{D}
x)dzdx \nonumber \\
 \ &\ & = C^{1}/2 \int_{K/ \delta}^{\sqrt{t/ M}}  \frac{1}{x}\exp(-\frac{\pi ^{2} t}{2 x^{2}})
\exp (-\lambda_{D} x)(\delta ^{2} x^{2}-K^{2})dx \nonumber \\
\ & \ & \geq \delta ^{2}C^{1}/2 \int_{  K/ \delta}^{\sqrt{t/ M}}x
\exp (-\frac{\pi ^{2} t}{2 x^{2}}) \exp (-\lambda_{D} x)dx - I,
\nonumber
\end{eqnarray}
where $$I =(C^{1}/2)K^{2}\int_{0}^{\infty} \frac{1}{x}\exp (-\frac{\pi
^{2} t}{2 x^{2}}) \exp (-\lambda_{D} x)dx .$$
 From the Laplace
transform method, equation (\ref{laplacemethod}), there exists
$C_{0}>0$ such that as $t\to\infty$,
\begin{equation}\label{lateknown}
I\sim
C_{0}t^{-1/6}\exp(-\frac{3}{2}\pi^{2/3}\lambda_{D}^{2/3}t^{1/3}).
\end{equation}
From equation (\ref{generalast11}) as $t\to \infty$,
\begin{equation}\label{exactast}
\int_{ 0}^{\infty} x \exp (-\frac{\pi ^{2} t}{2 x^{2}}) \exp
(-\lambda_{D} x)dx \sim
2\sqrt{\frac{\pi}{3}}(\frac{\pi^{2}}{2})^{1/2}\lambda_{D}^{-1}
t^{1/2} \exp (-\frac{3}{2}\pi^{2/3}\lambda_{D}^{2/3}t^{1/3}).
\end{equation}
Now for some $c_{1}>0$,
\begin{eqnarray}
\ &\ &
\int_{ 0}^{ K/ \delta}  x\exp (-\frac{\pi ^{2}
t}{2x^{2}}-\lambda_{D} x)
 dx
\nonumber  \\
\ &\ & \leq  e^{- \pi^2 \delta ^{2}t /2K^{2}}\int_{ 0}^{ K/
\delta}x\exp (-\lambda_{D} x)dx
\lesssim  e^{-c_{1} t},\label{exactast1}
\end{eqnarray}
and
\begin{eqnarray}
\ & \ & \int_{ \sqrt{t/M}}^{\infty} x \exp (-\frac{\pi ^{2} t}{2
x^{2}})
\exp (-\lambda_{D} x)dx
 \leq  \int_{ \sqrt{t/M}}^{\infty}
x\exp (-\lambda_{D} x)dx  \nonumber \\
  \ &\ & = ( \sqrt{t/M} \lambda_{D}^{-1}+\lambda_{D}^{-2})\exp (-\lambda_{D}\sqrt{t/M} ). \label{exactast2}
\end{eqnarray}
Now from equations (\ref{lateknown})-(\ref{exactast2}) we get
\begin{equation}\label{exactlowerbound1}
\liminf_{t\to\infty}  t^{-1/2} \exp
(\frac{3}{2}\pi^{2/3}\lambda_{D}^{2/3}t^{1/3}) P_{z}[\tau
_{D}(Z)>t]\geq \delta ^{2}(C^{1}/2)
2\sqrt{\frac{\pi}{3}}(\frac{\pi^{2}}{2})^{1/2}\lambda_{D}^{-1}.
\end{equation}
For the upper bound for $P[\tau _{D}(Z)>t]$ from equation
(3.10) in \cite{deblassie},
\begin{equation}
P_{z}[\tau _{D}(Z) > t]=2\int_{0}^{\infty} \!\int_{u}^{\infty}
P_{\frac{u}{u+v}} [\eta _{(0,1)}>\frac{t}{(u+v^{2})}]
f(u)f(v)dvdu\ .
\end{equation}
 We define the following sets that make up the domain of integration,
\begin{eqnarray}
 & A_{1}& =\{(u,v):v\geq u\geq 0, \  u+v\geq \sqrt{t/M}
 \},\nonumber\\
 & A_{2} & =\{ (u,v): \ u\geq 0,\ v\geq K, \ u\leq v,  \  u+v\leq \sqrt{t/M}
 \},\nonumber\\
& A_{3} & =\{ (u,v):\  0 \leq u\leq v \leq K\}.\nonumber
\end{eqnarray}
Over the set $A_{1}$ we have for some $c>0$,
\begin{eqnarray}
\ & \ & \int\!\int_{A_{1}} P_{\frac{u}{u+v}} [\eta _{(0,1)}>\frac{t}{(u+v)^{2}}] f(u)f(v)dvdu\nonumber\\
\ &\ & \leq \int\!\int_{A_{1}}  f(u)f(v)dvdu\leq
\exp(-c\sqrt{t/M}).\label{upperA1}
\end{eqnarray}
The equation (\ref{upperA1}) follows from the distribution of
$\tau_{D}$ from Lemma 2.1 in \cite{nane}.

Since on $A_{3}$, $t/(u+v)^{2}\geq M$,
\begin{eqnarray}
\ & \ & \int\!\int_{A_{3}} P_{\frac{u}{u+v}} [\eta _{(0,1)}>\frac{t}{(u+v)^{2}}] f(u)f(v)dvdu\nonumber\\
\ &\ & \leq \int_{0}^{K}\!\int_{0}^{K} \exp (-\frac{\pi ^{2} t}{2 (u+v)^{2}}) f(u)f(v)dvdu.\nonumber\\
\ &\ & \leq
\exp(-\frac{\pi^{2}t}{8K^{2}})\int_{0}^{K}\!\int_{0}^{K}f(u)f(v)dvdu\leq
\exp(-\frac{\pi^{2}t}{8K^{2}}) .\label{upperA3}
\end{eqnarray}
Let $C_{1}=C_{1}(z)=2\pi (4/\pi)A(z)^{2}(1+\epsilon)^{3}$. For the
integral over $A_{2}$ we get,
\begin{eqnarray}
\ & \ & \int\!\int_{A_{2}} P_{\frac{u}{u+v}} [\eta _{(0,1)}>\frac{t}{(u+v)^{2}}] f(u)f(v)dvdu\nonumber\\
\ &\ \leq &  C_{1}\int_{0}^{K}\!\int_{K}^{\sqrt{t/M}-u} f(u)\exp (-\frac{\pi ^{2} t}{2 (u+v)^{2}}-\lambda_{D}v)dvdu\nonumber\\
 & + &  C_{1}\int_{K}^{1/2\sqrt{t/M}}\!\int_{u}^{\sqrt{t/M}-u}\frac{u}{u+v} \exp (-\frac{\pi ^{2} t}{2 (u+v)^{2}}-\lambda_{D}(u+v))dvdu \nonumber\\
\ &\ =& I+II. \label{upperA2}
\end{eqnarray}
Changing variables $u+v=z$, $u=w$
\begin{eqnarray}
\ &\ I &  = \int_{0}^{K}\!\int_{K}^{\sqrt{t/M}-u} \exp (-\frac{\pi
^{2} t}
{2 (u+v)^{2}}) f(u)\exp (-\lambda_{D}v)dvdu\nonumber\\
\ & \leq & \int_{0}^{K}\!\int_{w+K}^{\sqrt{t/M}} \exp (-\frac{\pi
^{2} t}{2 z^{2}}) f(w)\exp (-\lambda_{D}z)
\exp(\lambda_{D}w)dzdw\nonumber\\
\ & \leq & \exp(\lambda_{D}K)\int_{0}^{K}f(w)dw \int_{0}^{\infty} \exp (-\frac{\pi ^{2} t}{2 z^{2}})
\exp (-\lambda_{D}z)dz\nonumber\\
\ & \lesssim \ & t^{1/6} \exp
(-\frac{3}{2}\pi^{2/3}\lambda_{D}^{2/3}t^{1/3}).\label{upperA21}
\end{eqnarray}
Equation (\ref{upperA21}) follows from equation (\ref{newasymptotic}),
with $a=\pi^{2}/2$, $b=\lambda_{D}$.

Changing variables $u+v=z$, $u=w$
\begin{eqnarray}
\ &\ II &  = C_{1}\int_{K}^{1/2\sqrt{t/M}}\!\int_{u}^{\sqrt{t/M}-u}
 \frac{u}{(u+v)}\exp (-\frac{\pi ^{2} t}{2 (u+v)^{2}}-\lambda_{D}(u+v))dvdu \nonumber\\
\ & \ & \leq C_{1}\int_{K}^{1/2\sqrt{t/M}}\!\int_{2w}^{\sqrt{t/M}}
\frac{w}{z}\exp (-\frac{\pi ^{2} t}{2 z^{2}}-\lambda_{D}z)dzdw \nonumber \\
\ & \ & =C_{1}\int_{2K}^{\sqrt{t/M}}\!\int_{K}^{z/2}
\frac{w}{z}\exp (-\frac{\pi
^{2} t}{2 z^{2}}-\lambda_{D}z)dwdz \label{changevar}\\
\ & \ & \ \leq C_{1}/8\int_{2K}^{\sqrt{t/M}}z \exp (-\frac{\pi
^{2} t}{2 z^{2}}-\lambda_{D}z)dz\nonumber\\
 \ & \ & \ \leq
  (1+\epsilon)(C_{1}/8) 2\sqrt{\frac{\pi}{3}}(\frac{\pi^{2}}{2})^{1/2}\lambda_{D}^{-1}t^{1/2}(-\frac{3}{2}\pi^{2/3}\lambda_{D}^{2/3}t^{1/3}).
\label{upperA22}
\end{eqnarray}

Equation (\ref{changevar}) follows by changing the order of the
integration. And equation (\ref{upperA22}) follows from equation
(\ref{generalast11}).

Now from equations (\ref{upperA1}), (\ref{upperA3}),
(\ref{upperA21}) and (\ref{upperA22}) we obtain
\begin{equation}\label{exactupperbound}
\limsup_{t\to\infty}t^{-1/2}\exp(\frac{3}{2}\pi^{2/3}\lambda_{D}^{2/3}t^{1/3})P_{z}[\tau
_{D}(Z)>t]\leq (1+\epsilon)(\frac{C_{1}}{8})
2\sqrt{\frac{\pi}{3}}(\frac{\pi^{2}}{2})^{1/2}\lambda_{D}^{-1}.
\end{equation}
Finally, from equations (\ref{exactlowerbound1}) and
(\ref{exactupperbound}) and letting $\epsilon\to 0$, $\delta\to 1/2$ ,
\begin{eqnarray}
 2C(z)  & \leq &
\liminf_{t\to\infty} t^{-1/2}
\exp(\frac{3}{2}\pi^{2/3}\lambda_{D}^{2/3}t^{1/3})P_{z}[\tau_{D}(Z)>t]\nonumber\\
& \leq  & \limsup_{t\to\infty}
t^{-1/2}\exp(\frac{3}{2}\pi^{2/3}\lambda_{D}^{2/3}t^{1/3})
P_{z}[\tau_{D}(Z)>t]\leq  \pi C(z),\nonumber
\end{eqnarray}
where $ C(z)=\lambda_{D}\sqrt{2\pi /3}\left(
\psi(z)\int_{D}\psi(y)dy\right) ^{2}. $
\end{proof}
\section{ The process $Z_{t}^{1}$; Brownian-time Brownian motion}\label{BTBM}
In this section we study Brownian-time Brownian motion (BTBM),
$Z_{t}^{1}$ started at $z \in \RR{R}$.
 Let $X_{t}$ and $Y_{t}$
 be two  independent
one-dimensional Brownian motions, all started at $0$. BTBM is
defined to be $Z_{t}^{1}=x+X(|Y_{t}|)$. In $\RR{R}^{n}$, we
require $X$ to be independent one dimensional iterated Brownian
motions. If $D\subset \RR{R} ^{n}$  is an open set, write
$$
\tau_{D}(z)=\inf \{ t\geq 0: \ X_{t}+z\notin D\},
$$
and if $I\subset \RR{R}$ is an open interval, we write
$$
\gamma_{I}=\inf \{ t\geq 0: \ |Y_{t}|\notin I\},
$$
and
$$
\eta_{I}=\inf \{ t\geq 0: \ Y_{t}\notin I\}.
$$
Let $\tau_{D}(Z^{1})$ stand for the first exit time of BTBM from $D$. We have by the continuity of paths
\begin{equation}\label{probeqisop}
P_{z}[\tau _{D}(Z^{1}) > t]=P[\eta(-\tau _{D}(z) , \tau _{D}(z) ) >t]. 
\end{equation}

\begin{theorem}\label{asymptotictheorem}
Let $0<\beta $. Let $\xi$ be a positive random variable such that
$$-\log P[\xi > t] \sim c t^{\beta}, \ as\  t\rightarrow \infty .$$ If
$\xi $ is independent of the Brownian motion $Y,$ then
\[
-\log P[\eta_{(-\xi,\xi)} >t] \sim
2^{-\frac{2\beta}{2+\beta}}(\frac{2+\beta}{2})c^{2/ (2+\beta)}
\beta^{-  \beta / (2+ \beta) } \pi ^{2 \beta/ (2+\beta)} t^{\beta
/ (2+ \beta)},
\]
as $t \rightarrow  \infty$.
\end{theorem}

\begin{proof}
The proof follows similar to the proof of Theorem 3.1 in Nane
\cite{nane},by integration by parts,
\begin{equation}\label{Z1representation}
P[\eta_{(-\xi,\xi)} >t]=\int_{0}^{\infty}  \frac{d }{d u} P_{0}
(\eta _{(-u,u)}>t)   P[\xi>u]du.
\end{equation}
We use the distribution of $\eta _{(-u,u)}$ given in (\ref{distr1}).
We use the
asymptotics from equation (\ref{lowerdiffisop}) on the set $A=\{
u>0:\ \ K\leq u \leq \sqrt{t/M}\}$. For the lower bound we use
Lemma 2.4 in \cite{nane}, but for the upper bound we use deBruijn Tauberian Theorem as in \cite[Lemma 2.2]{nane}.
\end{proof}

From Theorem \ref{asymptotictheorem} we obtain similar results for
the asymptotic distribution of the first exit time of $Z^{1}$ from
the interior of several open sets $D\subset \RR{R}^{n}$.

\begin{corollary}
Let $ 0 < \alpha <1$. Let $P_{\alpha}=\{ (x,Y)\in \RR{R} \times
\RR{R}^{n-1}:
 x>0, |Y|<Ax^{\alpha} \}.
$ Then for $z \in P_{\alpha}$,
 $$\lim_{t\to\infty}t^{-(\frac{ 1- \alpha }{3+\alpha})}\log P_{z}[\tau _{\alpha}(Z^{1})> t]
 =-2^{(\frac{2\alpha -2}{3+\alpha})}(\frac{3+\alpha }{2+ 2 \alpha } )
(\frac{1+\alpha }{1-\alpha})^{(\frac{1- \alpha }{3+ \alpha} )} \pi
^{(\frac{2-2\alpha }{3+\alpha } )} l^{  (\frac{2+2\alpha
}{3+\alpha }) } ,
$$
where $l$ is the limit given by \ref{brownianlimit}.
\end{corollary}

\begin{corollary}
Let $D\subset \RR{R}^{2}$ be a twisted domain with growth
 radius $\gamma r^{p}$, $\gamma >0,$ $0<p<1$.
Then for $z\in D$,
$$ \lim_{t\to\infty} t^{-(\frac{1-p}{p+3})}\log
P_{z}[\tau_{D}(Z^{1})>t]= -2^{(\frac{2p -2}{3+p})}(\frac{3+p }{2+
2 p } ) (\frac{1+p}{1-p})^{(\frac{1-p}{3+p})} \pi
^{(\frac{2-2p}{3+p } )} l_{1}^{  (\frac{2+2p }{3+p }) } ,
$$
where $l_{1}$ is the limit given by the limit in \cite[Theorem
1.1]{DSmits}.
\end{corollary}

\begin{remark}
Notice that there is only a constant difference in the limit of
the asymptotic distribution of $\tau_{D}(Z)$ and that of
$\tau_{D}(Z^{1})$ (compare with the results in Nane \cite{nane}.)
\end{remark}
\begin{proof}[Proof of Theorem \ref{boundeddomain2}]

From equations  (\ref{distr1}), (\ref{eigenvalue-ex}) and  (\ref{probeqisop})
\begin{eqnarray}
& \ & P_{z}[\tau_{D}(Z^{1})\leq t]=\int_{0}^{\infty}  P_{0}[\eta _{(-u,u)}>t] f(u)du\label{asymptoticz1}\\
&  = & \frac{4}{\pi }\sum_{k=1}^{\infty}\sum_{n=0}^{\infty}
\frac{(-1)^{n}}{2n+1}\lambda_{k}\psi_{k}(z)\int_{D}\psi_{k}(y)dy
\int_{0}^{\infty}\exp (-\frac{(2n+1)^{2}\pi ^{2}t}{8 u^{2}} -
\lambda_{k} u )  du.\nonumber
\end{eqnarray}
From equation (\ref{newasymptotic}), for each $n, \ k$ we have with
$a=\frac{(2n+1)^{2}\pi ^{2}}{8 }$ and $b=\lambda_{k}$
\begin{eqnarray}
 & \ & \int_{0}^{\infty}\exp (-\frac{(2n+1)^{2}\pi ^{2}t}{8 u^{2}}  - \lambda_{k} u )  du\nonumber\\
&  \sim  &
 \pi^{5/6}2^{1/6}3^{-1/2}(2n+1)^{1/3}\lambda_{k}^{-2/3}t^{1/6}\exp(-\frac{3}{2}(2n+1)^{2/3}\pi^{2/3}\lambda_{k}^{2/3}2^{-2/3}t^{1/3}).\nonumber
\end{eqnarray}
With this, equation (\ref{asymptoticz1}) becomes
\begin{eqnarray}
 & \ & \int_{0}^{\infty}  P_{0}[\eta _{(-v,v)}>t] f(v)dv  \label{asymptoticforz1}\\
&  \sim  & \frac{4}{\pi }\sum_{k=1}^{\infty}\sum_{n=0}^{\infty}
\frac{(-1)^{n}}{2n+1}
\lambda_{k}\psi_{k}(z)\int_{D}\psi_{k}(y)dy \nonumber\\
&  \times &
\pi^{5/6}2^{1/6}3^{-1/2}(2n+1)^{1/3}\lambda_{k}^{-2/3}t^{1/6}\exp(-\frac{3}{2}2^{-2/3}
(2n+1)^{2/3}\pi^{2/3}\lambda_{k}^{2/3}t^{1/3}).\nonumber
\end{eqnarray}
To get the desired result we must prove that the following series
converge absolutely, which implies that the first term in the
series in (\ref{asymptoticforz1}) is the dominant term,
$$
\sum_{k=1}^{\infty}\sum_{n=0}^{\infty}
 (2n+1)^{-2/3}\lambda_{k}^{1/3}\exp(-\frac{3}{2}2^{-2/3}(2n+1)^{2/3}\pi^{2/3}\lambda_{k}^{2/3}\delta/2)<\infty.
$$
The series in $n$ for $k$ fixed
$$ \sum_{n=0}^{\infty} (2n+1)^{-2/3}
 \exp(-\frac{3}{2}2^{-2/3}(2n+1)^{2/3}\pi^{2/3}\lambda_{k}^{2/3}\delta/2) $$
$$\leq  \frac{\exp(-\frac{3}{2}2^{-2/3}\pi^{2/3}\lambda_{k}^{2/3}\delta/2)}{1-\exp (-\frac{3}{2}\pi^{2/3}\lambda_{1}^{2/3}\delta/2)}.$$
Since for $\delta >0$,
$$
\sum_{k=1}^{\infty}
 \exp(-\frac{3}{2}2^{-2/3}\pi^{2/3}\lambda_{k}^{2/3}\delta/3) \leq \infty ,
$$
we are done. This follows from the Weyl's asymptotic formula for
the eigenvalues $\lambda_{k}$, $\lambda_{k}\geq C_{n,D} k^{n/2}$,
see I. Chavel \cite{chavel}, where $C_{n,D}$ depends only the
dimension $n$, and the domain $D$, independent of $k$. From above
equation (\ref{asymptoticforz1}) the constant
$C(\lambda_{D})=\pi^{-1/6}2^{13/6}3^{-1/2}\lambda_{D}^{1/3}$,
where $\lambda_{D}=\lambda_{1}$ is the first eigenvalue of the
Dirichlet Laplacian in $D$.
\end{proof}
\section{Comparison of  IBM and BTBM}

\begin{proof}[ Proof of Theorem \ref{IBM1IBM2}]

From equation (3.10) in \cite{deblassie} we get
\begin{eqnarray}
 P_{z}[\tau _{D}(Z) > t]
\ & = & 2\int_{0}^{\infty} \!\int_{u}^{\infty} P_{0}[\eta _{(-u,v)}>t] f(u)f(v)dvdu \nonumber \\
\ & \leq & 2 \int_{0}^{\infty} \!\int_{u}^{\infty} P_{0}[\eta _{(-v,v)}>t] f(u)f(v)dvdu \label{goodbound} \\
\ & \leq & 2 \int_{0}^{\infty} \!\int_{0}^{\infty} P_{0}[\eta _{(-v,v)}>t] f(u)f(v)dvdu\nonumber\\
\ & = & 2 \int_{0}^{\infty}  P_{0}[\eta _{(-v,v)}>t] f(v)dv\nonumber\\
\ & = & 2P_{z}[\tau _{D}(Z^{1}) > t] .\label{goodbound2}
\end{eqnarray}
The inequality (\ref{goodbound}) follows from the fact that $(-u,v)\subset (-v,v)$. The
equality (\ref{goodbound2}) follows from equation (\ref{probeqisop}).

\end{proof}

Let $\phi$ be an increasing function. If we multiply the
inequality in the Theorem \ref{IBM1IBM2} by the derivative of $ \phi$ and
integrate in time we get
$$
E_{z}(\phi (\tau_{D}(Z)))\leq  2 E_{z}(\phi (\tau_{D}(Z^{1}))).
$$
In particular, for $p\geq 1$,
$$
E_{z}((\tau_{D}(Z))^{p})\leq  2 E_{z}((\tau_{D}(Z^{1}))^{p}).
$$

\section{Asymptotics}\label{approximation}

In this Section we will prove some lemmas that were used in
section 3 and section 4. The following lemma is proved in
\cite[Lemma A1]{deblassie} (it also follows from more general
results on
``intrinsic ultracontractivity").  We include it for completeness.

\begin{lemma}\label{lemmaA.1}
As $t\rightarrow \infty$,
\[
P_{x}[\eta _{(0,1)} >t]\  \sim \  \frac{4}{\pi}e^{-\frac{\pi ^{2}
t}{2} } \sin \pi x, \ \  \mathrm{uniformly} \ \ \mathrm{for}\ \
x\in (0,1).
\]
\end{lemma}

We  will next prove the similar results to Nane \cite[Lemma
4.2]{nane} that will be used for the process $Z^{1}$.
\begin{lemma}\label{isoplowasymptotic} Let $B=\{ u>0: \ t/ u^{2} >M\}$ for $M$ large.
 Then on $B$,
\begin{equation}
 \frac{d }{du} P_{0}[\eta _{(-u,u)} > t]\sim \
\exp(-\frac{\pi ^{2} t}{8 u^{2}}) \frac{\pi t}{u^{3}}.
\label{lowerdiffisop}
\end{equation}
\end{lemma}

\begin{proof}
If we differentiate
$
P_{0}[\eta _{(-u,u)} >t]$ which is given in (\ref{distr1})
we get
$$
\frac{d}{du}P_{0}[\eta _{(-u,u)} >t]=\frac{\pi
t}{u^{3}}\sum_{n=0}^{\infty}(2n+1)(-1)^{n}\exp
(-\frac{(2n+1)^{2}\pi ^{2}}{8 u^{2}} t).
$$
The result follows from this.
\end{proof}

\textbf{Acknowledgments.} I would like to thank  Professor Rodrigo
 Ba\~{n}uelos, my academic advisor, for suggesting this problem to me and for  his guidance on this paper. I also would like to thank anonymous referee for his
 help with the presentation of the paper.

\end{document}